\documentclass[12pt]{article}

\usepackage[intlimits,tbtags]{amsmath}
\usepackage{amsfonts}
\usepackage{amsthm}

\usepackage{amssymb}
\usepackage{cite}
\usepackage{latexsym}
\usepackage[cmtip,arrow]{xy}
\usepackage{pb-diagram,pb-xy}

\usepackage{pst-node}
\usepackage{pstricks-add}
\usepackage{pstricks}

\newcommand{\D}{{d}}
\newcommand{\I}{{i}}
\newcommand{\E}{{e}}
\newcommand{\bbZ}{{\mathbb{Z}}}
\newcommand{\bbC}{{\mathbb{C}}}
\newcommand{\bbR}{{\mathbb{R}}}

\newcommand{\bbT}{{\mathbb{T}}}
\newcommand{\calA}{{\mathcal{A}}}

\newcommand{\calX}{{\mathcal{X}}}

\newcommand{\calU}{{\mathcal{U}}}

\DeclareMathOperator{\ev}{ev}
\DeclareMathOperator{\cv}{cv}
\DeclareMathOperator{\inv}{inv}
\DeclareMathOperator{\diag}{diag}
\DeclareMathOperator{\graph}{graph}
\DeclareMathOperator{\Hom}{Hom}
\DeclareMathOperator{\id}{id}
\DeclareMathOperator{\Id}{Id}

\newcommand{\op}{{\mathrm{op}}}

\DeclareMathOperator{\lcm}{lcm}

\DeclareMathOperator{\pr}{pr}

\DeclareMathOperator{\Cois}{Cois}
\DeclareMathOperator{\Lag}{Lag}
\DeclareMathOperator{\Iso}{Iso}

\newcommand{\Cat}[1]{{\mathsf{#1}}} 

\newcommand{\LieGrpdPrBibu}{\Cat{LieGrpdPrBibu}}

\newcommand{\Bimo}[1]{{\boldsymbol{#1}}} 

\newcommand{\Abas}{{a}}

\newcommand{\Bimod}[1]{{\boldsymbol #1}} 
\DeclareMathOperator{\Motimes}{{\otimes_{\Bimod{\Delta}}}} 

\newtheorem{Theorem}{Theorem}

\newtheorem{Definition}{Definition}

\theoremstyle{definition}

\makeatletter

\@namedef{dgo@||}{\let\dg@VECTOR=\dg@doublevvector\relax%
 \multiply\dg@LBLOFF 3\relax \divide\dg@LBLOFF 2}%
\def\dg@doublevvector(#1,#2)#3{%
  \begingroup
  \dg@XTEMP = \dg@LBLOFF\relax
  \divide\dg@XTEMP by 2
  \begin{picture}(0,0)%
     \thinlines
     \put(\dg@XTEMP,0){\vector(#1,#2){#3}}%
     \put(-\dg@XTEMP,0){\vector(#1,#2){#3}}%
  \end{picture}%
  \endgroup}

\@namedef{dgo@--}{\let\dg@VECTOR=\dg@doublehvector\relax%
 \multiply\dg@LBLOFF 3\relax \divide\dg@LBLOFF 2}%
\def\dg@doublehvector(#1,#2)#3{%
  \begingroup
  \dg@XTEMP = \dg@LBLOFF\relax
  \divide\dg@XTEMP by 2
  \begin{picture}(0,0)%
     \thinlines
     \put(0,\dg@XTEMP){\vector(#1,#2){#3}}%
     \put(0,-\dg@XTEMP){\vector(#1,#2){#3}}%
  \end{picture}%
  \endgroup}

\makeatother

\psset{unit=0.6}
\psset{linewidth=0.5pt}
\psset{fillcolor=lightgray,fillstyle=solid,linecolor=black}
\newgray{insidegray}{0.9}

\def\Disk(#1,#2){
 \psellipse(! #1 #2 )(2,0.5)
}

\def\Cone(#1,#2){
 \psline{-}(! #1 #2 )(! #1 -1.5 add #2 -1.5 add)(! #1 +1.5 add #2 -1.5 add)(! #1 #2 )
 \psellipse[fillcolor=insidegray](! #1 #2 -1.5 add)(1.5,0.25)
}

\def\ConeB(#1,#2){
 \psline{-}(! #1 #2 )(! #1 -1 add #2 -2 add)(! #1 +1 add #2 -2 add)(! #1 #2 )
 \psellipse[fillcolor=insidegray](! #1 #2 -2 add)(1.0,0.25)
}

\def\Ring(#1,#2){
 \psellipse(! #1 #2 )(2,0.5)
 \psellipse[fillcolor=white](! #1 #2 )(! 2 0.8 mul 0.5 0.8 mul)
}

\def\Orna(#1,#2){
 \pnode(! #1 #2  +2 add ){Tip1}
 \pnode(! #1 #2  -2 add ){Tip2}
 \ncarc[arcangleA=55,arcangleB=30]{-}{Tip1}{Tip2}
 \ncarc[arcangleA=-55,arcangleB=-30]{-}{Tip1}{Tip2}
}

\long\def\symbolfootnote[#1]#2{\begingroup%
\def\thefootnote{\fnsymbol{footnote}}\footnote[#1]{#2}\endgroup}

\begin{document}

\vspace{3em}
\begin{center}

{\Large{\bf Group-like objects in Poisson geometry and algebra}}

\vspace{3em}

\textbf{Christian Blohmann$^{1,2}$ and Alan Weinstein$^{1}$}

\vspace{1em}

${}^1$Department of Mathematics, University of California\\
Berkeley, California 94720-3840\\[1em]

${}^2$International University Bremen, School of Engineering and Science, Campus Ring 1, 28759 Bremen, Germany
\\[1em]

\end{center}

\vspace{1em}

\symbolfootnote[0]{These notes are loosely based on the three lectures given by
Weinstein at the School on Poisson Geometry and Related Topics, Keio University, May
31--June 2, 2006.  We would like to thank Nathan George for the use of
his preliminary notes.}

\begin{abstract}
A group, defined as set with associative multiplication and inverse, is a natural structure describing the symmetry of a space. The concept of group generalizes to group objects internal to other categories than sets. But there are yet more general objects that can still be thought of as groups in many ways, such as quantum groups. We explain some of the generalizations of groups which arise in Poisson geometry and quantization: the germ of a topological group, Poisson Lie groups, rigid monoidal structures on symplectic realizations, groupoids, 2-groups, stacky Lie groups, and hopfish algebras.
\end{abstract}

\section{Introduction}

Every mathematician learns that a group is a set with an associative
multiplication admitting an identity and inverses.  But there are
other objects, such as group germs, Lie groups, Poisson groups, and
quantum groups, which qualify as groups in many senses, but are either
more or less than simply sets with operations satisfying the group
actions.
These notes will describe some of these objects, with an emphasis on
the role of groupoids, both as examples of group-like objects and as a
tool for describing the objects themselves.  A goal toward which we
strive (but which we will not reach) is to give a unified
``categorical'' notion which encompasses all of our examples.

Some of our examples will be explicitly geometric.  Others will be
algebraic, but may still be considered as geometric from the viewpoint
which identifies geometric objects with suitable algebras of functions
on them and, more generally, considers algebras, even noncommutative
ones, as if they were the functions on a space.  Even more abstract is
the view of spaces as represented by categories, such as the category
of representations of an algebra, or
the derived category of coherent sheaves on an algebraic variety.

\section{Symmetry groups}

The (global) symmetry of a space is described by a set of
transformations closed under composition and inversion.  Conversely,
Cayley's theorem states that every abstract group is also a group of
transformations, acting on itself by left multiplication.

At first, a group was just a set, and its structure morphisms were maps
between sets. But it is useful to consider groups which themselves
have some additional structure, for example a differentiable
structure. The product, unit, and inversion
 morphisms of the group are required to
respect this structure, in which case they have to be smooth maps. This
leads us naturally to the concept of a Lie group, which is in turn
an example of a notion of group internal to a category, meaning that the
group is given by an object and group structure morphisms in that
category. In this sense, an ordinary group is a group object in the
category of sets, a topological group in the category of topological
spaces, a Lie group in the category of manifolds, an algebraic group
in the category of algebraic varieties, and so on. This notion
works well for many categories in which the objects are spaces with
geometric or topological structure.

On the other hand, in categories of spaces with algebraic structure,
the group objects often turn out to be surprisingly rare. For example,
the group objects in the category of vector spaces are the vector
spaces themselves with the underlying abelian group structure, while
in the category of groups they are the abelian groups. In the category
of rings the only group object is the trivial ring with one element.

But there are other kinds of objects which we can naturally associate
to symmetries. For example, to the action of a finite group $G$ on a
finite set $S$, the Gelfand ``algebraization'' functor
associates the commutative algebras of functions $\calA(G)$ and
$\calA(S)$ to both the group and the space. This is a contravariant
functor, so the group product, unit, and inverse in the group become
the comultiplication, counit, and coinverse (antipode) of a Hopf
structure on $\calA(G)$. The action of the group on the set becomes a
coaction of the Hopf algebra $\calA(G)$ on the algebra $\calA(S)$.
Since we can recover the group $G$ from $\calA(G)$ as the set of
group-like elements and the space $S$ from $\calA(S)$ as the set of
characters, the description of the symmetry of $S$ in terms of the
Hopf algebra $\calA(G)$ is completely equivalent to the description in
terms of the group.  This is why such Hopf algebras, even
noncommutative ones, have been termed quantum groups. But quantum
groups are not group objects, at least not in any of the underlying
categories of vector spaces, of algebras, or of coalgebras.  For
instance, the coproduct is a map from  $\calA(G)$ to $\calA(G)\otimes \calA(G)$, but the
tensor product is not a product in the category of algebras (nor in
the dual, with arrows reversed).  Our conclusion is that there are yet
more general structures which we may associate with the concept of a
symmetry group.

In fact, there is an ample collection of structures which are
considered to be ``groups'' in the sense that they encode symmetries:
groupoids, inverse semi-groups, hypergroups, $n$-groups, Lie algebras,
Hopf algebras, etc. --- just to name a few. Our goal in these lectures is to
identify some of such group-like structures that arise in Poisson
geometry and quantization.

\section{Group objects}

A group germ is an example of a group object in a category which is
not (at least in its usual presentation)
a subcategory of the category of sets.

We define a {\bf topological germ} to be the collection of all pointed
topological spaces $(X,x)$ modulo the equivalence relation in which
two spaces $(X,x)$ and $(Y,y)$ are identified if $x=y$ and if this
common point has open neighborhoods in $X$ and $Y$ which are equal as
sets and homeomorphic with the induced topologies.  A morphism between
topological germs is an equivalence class of continuous maps between
representatives of the germs, where two maps are considered equivalent
if they agree on some neighborhood of the basepoint.  This category
admits products defined as products of representatives, and a terminal
object $1$ consisting of a single point.  In any such category, a {\bf
group object} is defined to be an object $\calU$
together with the structure morphisms of multiplication $m: \calU
\times \calU \rightarrow \calU$, unit $e: 1 \rightarrow \calU$, and
inverse $\inv: \calU \rightarrow \calU$, such that the following
diagrams  are commutative. The diagrams for associativity and the unit
are
\begin{equation}
\label{eq:AssocUnital}
\begin{diagram}
\node{\calU \times \calU \times \calU}
 \arrow{e,t}{m \times \id}
 \arrow{s,l}{\id \times m}
\node{\calU \times \calU}
 \arrow{s,r}{m}
\\
\node{\calU \times \calU}
 \arrow{e,t}{m}
\node{\calU}
\end{diagram}
\qquad
\begin{diagram}
\node{1 \times \calU}
 \arrow{s,l}{e \times \id}
\node{\calU}
 \arrow{w,t}{\cong}\arrow{e,t}{\cong}
 \arrow{s,r}{\id}
\node{\calU \times 1}
 \arrow{s,r}{\id \times e}
\\
\node{\calU \times \calU}
 \arrow{e,t}{m}
\node{\calU}
\node{\calU \times \calU}
 \arrow{w,t}{m}
\end{diagram}
\end{equation}
while that for the inverse is
\begin{equation}
\label{eq:Inverse}
\begin{diagram}
\node[3]{\calU}
 \arrow[1]{sww,t}{\diag}
 \arrow[1]{see,t}{\diag}
 \arrow[2]{s,r}{\eta}
\\[1]
\node{\calU \times \calU}
 \arrow[2]{s,l}{ \inv \times \id }
\node[4]{\calU \times \calU}
 \arrow[2]{s,r}{ \id \times \inv }
\\
\node[3]{1}
 \arrow[2]{s,r}{e}
\\
\node{\calU \times \calU}
 \arrow[1]{see,b}{m}
\node[4]{\calU \times \calU}
 \arrow[1]{sww,b}{m}
\\[1]
\node[3]{\calU}
\end{diagram}
\end{equation}
where $\diag$ is the diagonal morphism, that is, the unique morphism
which lifts the identity on $\calU$ to the product $\calU \times
\calU$:
\begin{equation}
\label{eq:ProductDiag}
\begin{diagram}
\node[2]{\calU \times \calU}
 \arrow{sw,t}{\pr_1}
 \arrow{se,t}{\pr_2}\\
\node{\calU} \node[2]{\calU} \\
\node[2]{\calU}
 \arrow{nw,b}{\id}
 \arrow{ne,b}{\id}
 \arrow[2]{n,lr}{\exists !}{\diag}
\end{diagram}
\end{equation}

Restricting the group structure of a Lie $G$ group to its
germ at the unit element $e$ yields such a
group object in the category of topological germs
 which is called the germ of the group at $g$. This
is no longer a group in the usual sense, because ``it has only one
point''.   In a similar way, one can define manifold germs and group
objects in the category thereof.  This is useful, for instance, when
one tries to integrate a Banach Lie algebra to a group.  The global
object does not always exist
\cite{do-la:espaces}\cite{es-ko:nonenlargible}, 
but its germ does (and is
unique up to isomorphism).  

\section{Poisson Lie groups}

What is a Poisson Lie group?   It is usually defined as a
Lie group $G$ with a Poisson structure such that the multiplication
morphism $G\times G\to G$ is a Poisson map.  It follows from this that
the unit map from a point to $G$ is a Poisson map, and the inversion
map is anti-Poisson.


If we try to define a Poisson Lie group as a group object in the
category of Poisson manifolds, a 
first problem is that the category of Poisson manifolds does not
admit categorical products.  The cartesian product $X \times Y$ does
not work, since Poisson maps $A\to X$ and $A\to Y$ yield a Poisson map
$A\to X \times Y$ only when the images in $C^\infty (A)$ of
$C^\infty(X)$ and $C^\infty(Y)$ Poisson commute.

But let's forget this for a moment and admit $X\times X$ as some kind
of product, if not a categorical one, and
assume that $X$ is a Poisson Lie semigroup, i.e. a
Poisson manifold with an associative
multiplication map $m: X \times X \rightarrow X$, which is a Poisson
map, and a Poisson unit
map, $e :1 \rightarrow X$.  If $X$ is a group, we know 
that inversion is an anti-Poisson map, i.e. a Poisson 
map $\inv: X^\op \rightarrow
X$.  The categorical defining property of inversion is commutativity
of the diagram~\eqref{eq:Inverse}, but we must put in the
opposite Poisson structure to get
\begin{equation}
\label{diagram-opinverse}
\begin{diagram}
\node[3]{X} \arrow[1]{sww,t}{\diag} \arrow[1]{see,t}{\diag}
 \arrow[2]{s,r}{\eta} \\[1] \node{X^\op \times X} \arrow[2]{s,l}{ \inv
 \times \id } \node[4]{X \times X^\op} \arrow[2]{s,r}{ \id \times \inv
 } \\ \node[3]{1} \arrow[2]{s,r}{e} \\ \node{X \times X}
 \arrow[1]{see,b}{m} \node[4]{X \times X} \arrow[1]{sww,b}{m} \\[1]
 \node[3]{X}
\end{diagram}
\end{equation}
But the diagonal map is not a Poisson map from $X$ to $X^\op \times X$
or $X \times X^\op$ unless the Poisson bivector is zero.  So the
inverse axiom does not have an evident interpretation in the Poisson category.

Another approach is to analyze the inversion map via
its graph, which is
\begin{equation}
 \graph(\inv)
 = \{(x,x^{-1}) \,|\, x \in X \}
 = \{ (x,y) \in X\times X \,|\, xy = e \} \,,
\end{equation}
i.e., the pull-back, $ (X \times X) \times_{m,X,e}
1$.  This can be expressed in terms of the graph of multiplication and
the opposite of the graph of the unit, $\graph(e)^\op = \{ (e,1) \}$,
which are coisotropic submanifolds,
\begin{equation}
 \graph(m) \in \Cois((X\times X) \times X^\op) \,,\quad
 \graph(e)^\op \in \Cois(X \times 1^\op) \,.
\end{equation}
The composition of these two coisotropic submanifolds, viewed as
Poisson relations, is again a coisotropic submanifold,
\begin{equation}
 \graph(\inv)
 = \graph(m) \times_X \graph(e)^\op \in \Cois(X \times X) \,,
\end{equation}
where we have used that $\Cois((X \times X) \times 1^\op) \cong
\Cois(X \times X)$.  

This makes it possible to define a Poisson Lie group in the following
way.  Starting with any Poisson semigroup, we may first require that $m$ be
transverse to $e$, so that the fiber product 
$
\graph(m) \times_X \graph(e)^\op \in \Cois(X \times X) 
$
is a coisotropic submanifold, and then require that the projection of
the manifold to one of the factors of $X\times X$ be a diffeomorphism.
This is still not completely ``categorical,'' but we will see below
that it has a useful algebraic analogue.

\section{Poisson Lie groups and symplectic realizations}

It is sometimes possible to describe a group structure on a mathematical
object as an extra structure on a category of  ``representations'' of
the object.

For a Poisson manifold $X$, one such category is that of the 
symplectic realizations, in which the objects are the Poisson maps from
symplectic manifolds to $X$ \cite{we:local}, and the morphisms are
symplectic maps forming commutative diagrams with these.  
We may think of these as
symplectic ``points'' of $X$, or as
geometric representations of $X$
in the sense that a Poisson map $J: S \rightarrow X$ induces a
representation of the Poisson Lie algebra of functions on $X$ by
that of $S$, or by
hamiltonian vector fields:
\begin{equation}
 C^\infty(X) \longrightarrow \calX(S)
 \,,\qquad
 f \longmapsto - X_{J^*(f)} \,.
\end{equation}
There always exists a symplectic realization which is a surjective
submersion \cite{co-da-we:groupoides}\cite{ka:analogues}, for which 
(up to locally constant functions) this
representation is faithful. Therefore, the collection of symplectic
realizations encodes all the structural information of $X$.  

Given two symplectic realizations $J_1:S_1 \rightarrow X$ and $J_2:S_2
\rightarrow X$ we can use a Poisson Lie structure on $X$ to construct
a product realization $J_1 \otimes J_2: S_1 \otimes S_2 := S_1 \times S_2 :
\rightarrow X$ by
\begin{equation}
\begin{diagram}
\node[2]{S_1 \times S_2}
 \arrow{sw,t}{J_1 \times J_2}
 \arrow{se,t}{J_1 \otimes J_2}
\\
\node{X \times X}
 \arrow[2]{e,t}{m}
\node[2]{X}
\end{diagram}
\quad.
\end{equation}
This multiplication of symplectic realizations is associative because
the multiplication on $X$ is. Furthermore, viewing the terminal object
in the category of Poisson manifolds $1 = \{ \mathrm{pt} \}$ as
a zero-dimensional symplectic manifold, the unit $e:1 \rightarrow X$ can
also be viewed as symplectic realization. It is the identity for the
product of realizations $e \otimes J = J = J \otimes e$, where we
identify $1 \times S = S = S \times 1$. In this way, the monoidal
structure on $X$ naturally equips the category of symplectic
realizations with a monoidal structure.

What structure is induced on the category of symplectic realizations
by the inverse on $X$? From the analogous algebraic situation, we
might expect that the inverse leads to a \emph{rigid} monoidal
structure \cite{ul:tannakian}.
We can try to define a dual symplectic realization by
\begin{equation}
 J^\vee : S^\vee
 \equiv S^\op \stackrel{J}{\longrightarrow}
 X^\op \stackrel{\inv}{\longrightarrow} X \,.
\end{equation}
But in the category of symplectic realizations of $X$ the evaluation map
would have to make the diagram
\begin{equation}
\label{eq:eval}
\begin{diagram}
\node{S \otimes S^\vee}
 \arrow{se,b}{J \otimes J^\vee}
 \arrow[2]{e,t}{\ev}
\node[2]{1}
 \arrow{sw,b}{e}
\\
\node[2]{X}
\end{diagram}
\quad,
\end{equation}
commutative, which is not possible unless $J$ maps all of $S$ to a
single point in $X$. If we want to equip the category of symplectic
realizations with a rigid structure, we will need a more general notion
of morphism.

\section{Generalized morphisms}

In many categories, morphisms are given by set theoretic maps, but we
may allow relations instead of just maps, especially when maps of a
certain type are characterized by properties of their graphs.
For instance, a smooth map $f: A \rightarrow B$ between symplectic manifolds
is symplectic if and only if its graph is a lagrangian submanifold of
$A \times B^\op$. Moreover, for the graph of the composition of $f$
with another map $g: B \rightarrow C$ we have
\begin{equation}
\label{eq:GraphCompose}
 \graph(g \circ f) = \graph(f) \circ \graph(g) := \graph(f) \times_{B} \graph(g) \,.
\end{equation}
This suggests allowing arbitrary lagrangian submanifolds of products
as morphisms, 
i.e., defining
\begin{equation}
 \Hom(A,B) := \Lag(A \times B^\op) \,,
\end{equation}
the symplectic $A$-$B$ relations. Note that, for the morphisms in
$\Lag(A,B)$, there is no natural distinction between source and target,
as there is for maps. A lagrangian submanifold $L$ of $A \times B^\op$
is the same as a lagrangian submanifold of $(A \times B^\op)^\op =
A^\op \times B \cong B \times A^\op$. So $L$ can be equivalently
viewed as a morphism from $B$ to $A$. This is why we prefer to denote
the source and target maps of a category by $l$ and $r$, because
everyone agrees what is left and right.

Using this generalized notion of morphisms we return to the symplectic realizations of a Poisson Lie group $X$. A generalized morphism between two symplectic realizations $J_1: S_1 \rightarrow X$ and $J_2: S_2 \rightarrow X$ is given by a lagrangian submanifold $L \in \Lag(S_1 \times S_2^\op)$ such that the following diagram commutes:
\begin{equation}
\begin{diagram}
\node[2]{L}
  \arrow{sw,t}{\pr_1}
  \arrow{se,t}{\pr_2}
\\
\node{S_1}
  \arrow{se,b}{J_1}
\node[2]{S_2}
  \arrow{sw,b}{J_2}
\\
\node[2]{X}
\end{diagram}
\end{equation}
Now we can try again to find an evaluation morphism from $S \otimes S^\vee$ to $1$ as in Eq.~\eqref{eq:eval}. What we need is a Lagrangian submanifold $L_{\ev} \in \Lag(S \times S^\op) \cong \Lag((S \times S^\op) \times \{\mathrm{pt}\} )$ such that for all $(s,s') \in L_{\ev}$ we have $J(s) J(s')^{-1} = e$. The natural lagrangian submanifold satisfying these requirements is the diagonal $L_{\ev} := \Delta_S = \{(s,s) \,|\, s \in S \}$. The same reasoning leads us to define the coevaluation morphism also by $L_{\cv} := \Delta_S$. Moreover, we have the same evaluation and coevaluation morphisms for $S^\vee \otimes S$.  It is easy to see that the morphism of symplectic realizations
\begin{equation}
  S 
  \xrightarrow{\cong}
  S \otimes 1 
  \xrightarrow{\Id_S \otimes L_{\cv}} 
  S \otimes (S^\vee \otimes S)
  \xrightarrow{\cong} 
  (S \otimes S^\vee) \otimes S
  \xrightarrow{L_{\ev} \otimes \Id_S}
  1 \otimes S
  \xrightarrow{\cong}
  S
\end{equation} 
is the identity morphism, which is also given by the diagonal $\Id_S =\Delta_S$. Going in an analogous way from $S^\vee$ via $S^\vee \otimes S \otimes S^\vee$ to $S^\vee$ we obtain the identity on $S^\vee$ which is also given by the diagonal. We are tempted to conclude that the category of symplectic realizations and generalized morphisms is rigid monoidal. However, there is a catch:

Unfortunately, the symplectic relations are  not really the morphisms
of a category; when
the projections to $B$ of 
elements in $\Lag(A,B)$ and $\Lag(B,C)$  intersect badly,
their composition as defined in
Eq.~\eqref{eq:GraphCompose} is not a manifold.
To avoid this, we can define $\Hom(A,B)$ to be $A\times B^\op$ itself,
rather than the set of lagrangian submanifolds therein.  The price we
pay is that the composition operation 
\begin{multline}
 \Hom(A,B) \times \Hom(B,C) = (A \times B^\op) \times (B \times C^\op) \\
 \longrightarrow A \times C^\op = \Hom(A,C) 
\end{multline}
is not a mapping of sets, but a relation, namely the
lagrangian submanifold of 
\begin{equation}
 \bigl( (A \times B^\op) \times (B \times C^\op) \bigr)
 \times \bigl( A \times C^\op \bigr)^\op
\end{equation}
consisting of the product $\{(a,b,b,c,a,c)|a\in A, b\in B, c\in C\}$
of three diagonals.  The result is what we have called in 
Section 5.2 of \cite{bu-we:poisson} a ``symplectic category,'' i.e. a
category internal to the ``category'' of symplectic relations.

\section{Groupoids and stacks}

Even more general than a relation between the sets $X$ and $Y$ is a 
``multi-relation'', i.e. a map from a set $M$ to $X \times Y$, which
might not be injective.  The best theory of such generalized morphisms
comes about when $M$ is acted upon by groupoids over $X$ and
$Y$.\footnote{It may not always be necessary to bring in the groupoids from
  the beginning; see for instance the use of ``bisubmersions'' in 
\cite{an-sk:holonomy}.}  The result is the theory of stacks, which we
will describe in its smooth version
\cite{be:stack}\cite{be-xu:differentiable}\cite{ts-zh:stack}.  
(See \cite{mr:stability} for the topological case.) 

 Roughly
speaking, a stack is a device to describe a ``bad'' quotient. 
Here is a simple example. Let $\bbZ_2$ act on an open
disc of unit radius in $\bbR^2$ by reflection at the origin, $1 \cdot
(x,y) = (-x,-y)$. The $\bbZ_2$ action is not free, because the origin
is a fixed point. Taking the quotient amounts to cutting the disc
along, say, the positive $x$-axis and rolling it up until you have two
layers at every point with the exception of the origin.
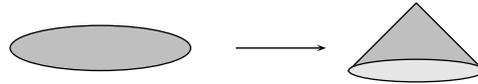
\begin{equation*}
\begin{pspicture}[0.5](-2,-1)(8,1)
 \Disk(0,0)
 \psline{->}(3,0)(5,0)
 \Cone(7,1)
\end{pspicture}
\end{equation*}
The result is the surface of a cone, which is no longer a
smooth manifold. You can do this for other finite cyclic groups, say,
$\bbZ_3$ where now $1 \in \bbZ_3$ acts by rotation by a third of the
full circle. Now you obtain a cone with a smaller opening angle.
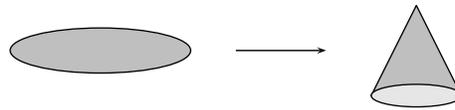
\begin{equation*}
\begin{pspicture}[0.5](-2,-1.5)(8,1)
 \Disk(0,0)
 \psline{->}(3,0)(5,0)
 \ConeB(7,1)
\end{pspicture}
\end{equation*}
 We can smoothly glue together two of these cones along an outer annulus
 of each disk to obtain a ``Christmas tree ornament'':
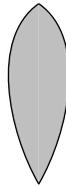
\begin{equation*}
\begin{pspicture}[0.5](-1,-2)(1,2)
 \Orna(0,0)
\end{pspicture}
\end{equation*}
This is an example for an orbifold, a manifold with lower dimensional
singularities which look locally like the quotient of $\bbR^n$ by a
finite group.  (The example above is still a manifold, topologically,
but this is not true for more general orbifolds.)
Like a manifold, an orbifold can be described by charts
with these quotient spaces as local models. A drawback of this
description is that there is no natural tangent bundle of an orbifold
which is itself an orbifold, and the definition of morphism is 
rather complicated and unnatural-looking.

A more effective way to describe a bad quotient space is to remember all
the gluings. This leads us to the concept of a groupoid. As a set, the
gluing groupoid for $D^2/\bbZ_2$ is $G_1 := \bbZ_2 \times
D^2$. There are two maps to the base $G_0 = D^2$, which map each
element $g$ of the groupoid to the two points which $g$ glues
together. We denote these maps by $l$ and $r$,
\begin{xalignat}{2}
 l\bigl(0,(x,y)\bigr) &= (x,y) \,,& r\bigl(0,(x,y)\bigr) &= (x,y) \,,\\
 l\bigl(1,(x,y)\bigr) &= (-x,-y) \,,& r\bigl(1,(x,y)\bigr) &= (x,y) \,,
\end{xalignat}
that is, $l(a,p) = a \cdot p$, $r(a,p) = p$ for $a \in \bbZ_2$ and $p
\in D^2$. We can compose two elements of the groupoid $(a,p)(a',p') =
(a+a',p')$ whenever $r(a,p) = p = a'\cdot p' = l(a',p')$.  This
construction can be extended to arbitrary group actions, but but there
are more general groupoids.  For example, the groupoid presenting the
Christmas ornament orbifold looks like this.
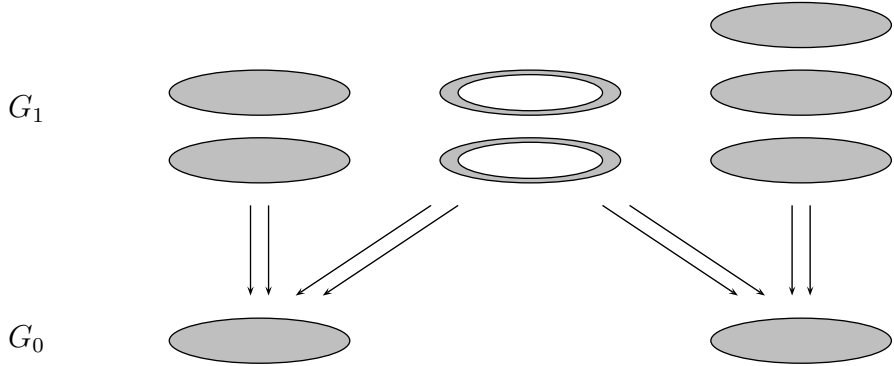
\begin{equation*}
\begin{matrix}
 \\ \\ G_1 \\ \\ \\ \\ \\ \\ G_0
\end{matrix}\qquad\qquad
\begin{pspicture}[0.5](-8,-5)(8,3.5)
 \Disk(-6,-4)
 \psline{->}(-6.2,-1)(-6.2,-3)
 \psline{->}(-5.8,-1)(-5.8,-3)
 \Disk(-6,0) \Disk(-6,1.5)
 \Ring(0,0)\Ring(0,1.5)
 \psline{->}(-2.2,-1)(-5.2,-3)
 \psline{->}(-1.6,-1)(-4.6,-3)
 \psline{->}(1.6,-1)(4.6,-3)
 \psline{->}(2.2,-1)(5.2,-3)
 \Disk(6,-4)\Disk(6,0)\Disk(6,1.5) \Disk(6,3)
 \psline{->}(5.8,-1)(5.8,-3)
 \psline{->}(6.2,-1)(6.2,-3)
\end{pspicture}
\end{equation*}
The $l$ and $r$ maps on the annuli are given by the embeddings into the disks.

All structure maps of the groupoid are smooth (and $l$ and $r$
submersions) so we have a Lie groupoid. In order to see what is
special about the tips of the Christmas ornament, we have to look at
the isotropy group $\Iso(p) := l^{-1}(p) \cap r^{-1}(p)$ for $p \in
G_0$, which for an action groupoid is the stabilizer of $p$. The
origins of the two disks of the base $G_0 = D^2 \cup D^2$ are the only
points which have non-trivial isotropy, $\bbZ_2$ for the left
disk and $\bbZ_3$ for the right disk. For orbifolds, the stabilizers
are by definition finite groups. For the groupoid it means that it is proper \'etale, i.e., the anchor map $(l,r): G_1 \rightarrow G_0
 \times G_0$, $g \mapsto (l(g),r(g))$ is proper \'etale. Indeed, we can give the following definition of an orbifold \cite{mo:orbifolds}:


\begin{Definition}
 An orbifold is a differentiable stack presented by a proper \'etale Lie groupoid.
\end{Definition}

Now we have to explain how we can associate a stack to a Lie groupoid. We can think of a Lie groupoid as a generalized equivalence relation
describing the quotient space of equivalence classes. As we have seen,
the actual quotient space $G_0/G_1$ of $G_1$-orbits in $G_0$ is
usually not a nice space, so we write $G_0/\!/G_1$ for the as-if
quotient the groupoid is thought to describe. The usual notion of
isomorphism of Lie groupoids is that of a diffeomorphism which is
compatible with all the structure maps. However, now two groupoids
should be considered to be equivalent if they present the same
quotient.

Given two groupoids $G$ and $H$ we must find a smooth way to associate
the $G_1$-orbits on $G_0$ with the $H_1$-orbits on $H_0$.  This
cannot be just a map $G_0/G_1 \rightarrow H_0/H_1,$ because these
quotient spaces are in general not manifolds, and so there is no notion of
smoothness. Our description of a bad quotient as a groupoid which is a
generalized equivalence relation suggests defining
morphisms between them in
an analogous manner. Now the two spaces we want to relate are the
bases $G_0$ and $H_0$ of the groupoids. A generalized relation
consists of triples $(x \stackrel{m}{\rightarrow} y)$ where we say
that $x \in G_0$ is related via $m$ to $y \in H_0$. We denote the set
of all such triples by $M$. The projections of the elements of $M$ on
the elements of the groupoid bases they relate, $l_M (x
\stackrel{m}{\rightarrow} y) := x$ and $r_M(x
\stackrel{m}{\rightarrow} y) := y$, are called the moment maps of
$M$. When no confusion can arise, we will drop the subscripts of the
moment maps.

We do not require $M$ to be a map from $G_0$ to $H_0$. For example, a
single pair
$x$ and $y$ can be related by several elements of $M$. But we want the
relation $M$ to descend to a map on the set of orbits.   For
notational reasons it is convenient to use left orbits in $G_0$ and
right orbits in $H_0$.  (Note that the
left orbits for a groupoids acting on its base
 are the same as the right orbits.) 
 For the generalized relation $M$ between elements $x \in
G_0$ and $y \in H_0$, to descend to a well-defined relation on the
orbits $G \cdot x$ and $y \cdot H$ we have to require that for all $g$
acting on $x$ and $h$ acting on $y$ we have elements $g\cdot m$ and
$m\cdot h$ of $M$ such that
\begin{equation}
 x \stackrel{m}{\rightarrow} y
 ~\Rightarrow~ (g \cdot x) \stackrel{g\cdot m}{\longrightarrow} y
 \qquad\text{and}\qquad
 x \stackrel{m}{\rightarrow} y
 ~\Rightarrow~ x \stackrel{m \cdot h}{\longrightarrow} (y \cdot h) \,.
\end{equation}
We want to be able to chose $g
\cdot m$ and $m \cdot h$ in a consistent way, such that we get maps $m
\mapsto g \cdot m$ and $m \mapsto m \cdot h$ which are compatible with
the groupoid structures, $g\cdot (g' \cdot m) = gg' \cdot m$ and $(m
\cdot h) \cdot h' = m \cdot hh'$, as well as, $(g \cdot m) \cdot h = g
\cdot (m \cdot h)$ whenever defined, i.e., we have two commuting
groupoid actions on $M$.  Such an object $M$ is called a groupoid
bibundle.

So far, $M$ only descends to a relation on $G_0/G_1 \times
H_0/H_.$ To obtain a function from $G_0/G_1$ to
$H_0/H_1$, we first have  to require that $l_M$ is
surjective, so that the function will be defined on all of
$G_0/G_1$. Second, for the relation to be a function, a given $x \in
G_0$ has to be related only to elements of one orbit in $H_0$;
\begin{equation}
 x \stackrel{m}{\rightarrow} y \quad\text{and}\quad x
 \stackrel{m'}{\rightarrow} y' \quad\Rightarrow\quad y \cdot h = y'
 \quad\Rightarrow\quad x \stackrel{m \cdot h}{\longrightarrow} y'
\end{equation}
for some groupoid element $h \in H_1$. Again,  we want to be able to
chose $h$ in a nice way, requiring that there is a unique $h$ such
that $m \cdot h = m'$.   This gives us a right principal
bibundle. Finally, we require all the structures to be smooth, that
is, $M$ is a manifold, the moment maps are smooth, and the groupoid
actions are smooth.

We can depict the situation by the diagram
\begin{equation}
\label{eq:smooth-nonsmooth}
\begin{diagram}
\node{G_1}
 \arrow{s,r,||}{r_G}
\node{M}
 \arrow{sw,b}{l_M}
 \arrow{se,b}{r_M}
\node{H_1}
 \arrow{s,l,||}{l_H}
\node{\parbox[r]{0ex}{\text{smooth}}}
\\
\node{G_0}
 \arrow{s}
\node[2]{H_0}
 \arrow{s}
\node{\parbox[r]{0ex}{\text{smooth}}}
\\
\node{G_0/G_1}
 \arrow[2]{e}
\node[2]{H_0/H_1}
\node{\parbox[r]{0ex}{\text{not smooth}}}
\end{diagram}
\end{equation}
Now let $G$, $H$, $K$ be Lie groupoids, $M$ a smooth right principal
$G$-$H$ bibundle and $N$ a smooth right principal $H$-$K$
bibundle. The composition of the induced functions on the quotient
spaces can be lifted to a smooth composition of the bibundles:
\begin{equation}
\begin{diagram}
\node[3]{M \circ N := (M \times_{H_0} N)/H}
 \arrow{sssww}\arrow{sssee}
\\
\node[3]{M \times_{H_0} N}
 \arrow{n}
 \arrow{sw}\arrow{se}
\\
\node{G_1}
 \arrow{s,||}{}
\node{M}
 \arrow{sw}
 \arrow{se,b}{r_M}
\node{H_1}
 \arrow{s,||}
\node{N}
 \arrow{sw,b}{l_N}
 \arrow{se}
\node{K_1}
 \arrow{s,||}
\\
\node{G_0}
\node[2]{H_0}
\node[2]{K_0}
\end{diagram}
\end{equation}
Here, the right $H$-action on $M \times_{H_0} N$ is given by $(m,n) \cdot
h = (m \cdot h, h^{-1} \cdot n)$. The $G$-action and the $K$-action
descend to actions on the quotient, since they both commute with the
$H$-action. The conclusion is that $M \circ N$ is a $G$-$K$ bibundle,
which we call the composition of $M$ and $N$. This composition is only
associative up to a biequivariant diffeomorphism of bibundles. 
This means that we should really be working in 
the weak 2-category having Lie groupoids as objects,
smooth right principal bibundles as 1-morphisms, and smooth
biequivariant maps of bibundles as 2-morphism. We denote this category
by $\LieGrpdPrBibu$.

A good way to study the generalized space described by a groupoid 
$G$ is the
Grothendieck approach of considering 
 all the morphisms from ordinary manifolds to
$G$, where  a manifold $X$ is described by
the groupoid $X_1 = X_0 = X$. The 
 collection of all such morphisms to $G$, more concretely described as

\begin{equation}
 \{ M \,|\, M \text{ a right principal $X$-$G$ bibundle},~ X \text{ a manifold} \} \,
\end{equation}
is denoted by $BG$
and called the classifying space of the groupoid $G$. 
This is a smooth stack in the usual
sense which is presented by the groupoid $G$. Two stacks $B G$ and $B
H$ are 
isomorphic if and only if the groupoids $G$ and $H$ are Morita
equivalent (Theorem 2.24 in \cite{be-xu:differentiable}). We thus get a 1-to-1 correspondence of isomorphism classes 
of presentable stacks and Morita equivalence classes of
groupoids. This is often stated as ``a stack is a groupoid up to
Morita equivalence''. But beware that the actual functor between the
category of stacks and the weak 2-category of groupoids and bibundles
is a weak 2-equivalence of 2-categories.

\section{Stacky Lie groups}

It can be shown that, in the weak 2-category $\LieGrpdPrBibu$ of Lie
groupoids and smooth right principal bibundles, all products exist and
the one-element groupoid $1$ is a terminal object.  (Note that, if the
bibundles are not required to be principal, this is no longer true.)
Since we are dealing with a weak 2-category, the categorical product
is associative only up to weak 1-isomorphisms, that is, up to Morita
equivalence of groupoids.

If we have products and a terminal objects we also have the notion of
group objects, which we will call stacky Lie groups
\cite{zh:ngroupoids,bl:stacky}. The question is whether in
$\LieGrpdPrBibu$ the notion of group objects is useful, as for
differentiable spaces, or uninteresting as for groups or
algebras. This is not easy to see because, on the one hand we think of
a groupoid as a generalized quotient space, but on the other hand a
Lie groupoid is itself an algebraic structure internal to the category
of manifolds. It turns out that while not many examples of truly
stacky groups have been studied until now, there are some particularly
interesting ones.

Consider the group $S^1$ and the dense subgroup which is the image of
the embedding $\bbZ \rightarrow S^1 \cong U(1)$, $k \mapsto e^{i
\lambda k}$ where, $\lambda/2\pi$ is an irrational number.  By abuse
of notation, we will also denote the subgroup itself by $\bbZ$.  The
quotient $S^1/\bbZ$ is an abelian group, in which we will denote the
multiplication by $m: S^1/\bbZ \times S^1/\bbZ \rightarrow
S^1/\bbZ$. Because the subgroup $\bbZ$ is dense, the quotient topology
is trivial. This suggests that we should work with the stack
$S^1/\!/\bbZ$ rather than with the actual quotient. We then try to
lift the multiplication map $m$ on the bottom level of the diagram
of the form \eqref{eq:smooth-nonsmooth} to the smooth top levels:
\begin{equation*}
\begin{diagram}
\node{(S^1 \times S^1)\times_{S^1/\bbZ} S^1}
 \arrow[2]{s}
 \arrow[2]{e}
 \arrow{se,t,..}{\exists!}
\node[2]{S^1}
 \arrow{s,A}
\\
\node[2]{(S^1/\bbZ \times S^1/\bbZ) \times_{S^1/\bbZ} S^1/\bbZ}
 \arrow{s}
 \arrow{e}
\node{S^1/\bbZ}
 \arrow{s,r}{\id}
\\
\node{S^1 \times S^1}
 \arrow{e,A}
\node{S^1/\bbZ \times S^1/\bbZ}
 \arrow{e,t}{m}
\node{S^1/\bbZ}
\end{diagram}
\end{equation*}
The inner pull-back $(S^1/\bbZ \times S^1/\bbZ) \times_{S^1/\bbZ} S^1/\bbZ$ is the graph of the group multiplication $m$, the pull-back projections being the range and image maps. Any object in the top left upper corner which makes the diagram commute can be viewed as a lift of the graph of $m$ to $S^1$, the pull-back being the universal lift. The diagonal arrow is the unique map which exists by the universal property of the inner pull-back. Explicitly, the pull-back is
\begin{equation*}
 (S^1 \times S^1)\times_{S^1/\bbZ} S^1 =
 \{(\theta_1,\theta_2, \theta)
 \in (S^1 \times S^1) \times S^1
 \,|\, \theta_1 + \theta_2 = \theta \bmod 2\pi\lambda \}
\end{equation*}
where $\theta$ is the angle representing $e^{i\theta} \in U(1) \cong
S^1$. This set can be identified with $S^1 \times S^1 \times \bbZ$ and
viewed as the graph of a multi-valued multiplication on $S^1$. It has
the structure of a smooth manifold and inherits smooth actions of the
action groupoid $G := S^1 \rtimes \bbZ \rightrightarrows S^1$
presenting the stack $S^1/\!/\bbZ$. We thus obtain the smooth right
principal $(G\times G)$-$G$ bibundle $E_m$ of multiplication. In an
analogous manner we construct the $1$-$G$ bibundle $E_e$ of the
identity element, and the $G$-$G$ bibundle of the inverse
$E_{\inv}$. It can be checked that the bibundles $E_m$, $E_e$, $E_{\inv}$
equip the groupoid $G$ with the structure of a stacky Lie group.

\section{Hopfish algebras}

The Gelfand Theorem tells us that a locally compact
topological space and its commutative algebra of continuous functions
vanishing at infinity contain the same structural information. In our
example, $S^1/\bbZ$, the quotient topology is trivial, so the only
continuous functions are constant.   But we are
working instead with the groupoid presenting the stack $S^1/\!/\bbZ$. What is
the algebra of continuous ``functions'' on $S^1/\!/\bbZ$? It is one of
the main ideas of noncommutative geometry that we should consider the
convolution algebra of the groupoid $G := S \rtimes \bbZ$ to be the
analogue of the algebra of functions
\cite{re:groupoid}\cite{co:ncgeometry}. For two compactly supported
functions $a$ and $b$ on $G$, the convolution product is
\begin{equation}
 (a * b)(g) = \int a(h) b( h^{-1} g) \D h \,,
\end{equation}
which looks just like the convolution algebra of a group. The
difference is that on a groupoid the product $h^{-1}g$ is only defined
if $l(h) = l(g)$, so for a given $g$ we have to integrate over the
left groupoid fiber of $r(g) = x$.  Since we have to do this for
all $g$, we need a whole family of Haar measures $\D h = \D_x
h$. Alternatively, we can work with half-densities instead of
functions, as in \cite{co:ncgeometry}.

Recalling that Morita equivalent groupoids have Morita equivalent
(but generally not isomorphic) convolution algebras, we also conclude
that Morita equivalent algebras should be thought of as representing
the ``same noncommutative space''.\footnote{One must use this
  identification with some care.  For instance, there are plenty of
  examples of nonisomorphic groups with isomorphic group algebras
  (e.g. pairs of finite abelian groups with the same number of
  elements); thinking of these groups as representing stacks of the
  form $BG$, one finds different stacks with the same ``algebra of
  functions''.  Even more different, it seems to us, 
 are the stacks presented by the group $Z_2$ and the trivial groupoid
 over a set with two elements, yet again they have isomorphic group algebras
 These examples suggest that it takes more than an algebra to
  make a ``noncommutative space,'' but it is not clear to us exactly
  what that ``more'' should be.}  

In our example, the Lie groupoid $G = S^1 \rtimes \bbZ$ is \'etale,
so the Haar measures are merely counting measures.  Thus,  
\begin{equation}
\label{eq:Convolution2}
 (a * b)(\theta,k)
 := \sum_{k' \in \bbZ} a(\theta + \lambda k',k-k') b(\theta,k') \,,
\end{equation}
for all compactly supported functions $a$ and $b$ on $S^1 \times
\bbZ$.   Among such functions are the standard Fourier basis
functions on the circle times the characteristic functions
of integers, $ a_{nl}(k, \theta)
 := \E^{\I n \theta} \delta_{lk},$
for which we obtain
\begin{equation}
 \Abas_{n_1 l_1} * \Abas_{n_2 l_2}
 = \E^{\I \lambda n_1 l_2} \Abas_{n_1 + n_2, l_1 + l_2} \,.
\end{equation}
Let us denote the algebra of finite linear combinations of these
 functions by $\calA$.  The closure of $\calA$ with respect to a
 suitable norm is known 
as the algebra continuous functions on a noncommutative
 torus,\footnote{By abuse of terminology, the algebra itself is often
 known as a noncommutative torus.}   which is
a deformation of the usual algebra of continuous functions on the 2-torus
\cite{ri:quantization}.

Now the quotient space $S^1/\bbZ$ is also a quotient
 group.  Since $\calA$ plays the role of 
an algebra of functions on this quotient,  we might expect
the group structure to translate into a Hopf structure on $\calA$. But
this is not the case. In fact, $\calA$ is simple, so
it does not even possess a counit. So what happened to the group
structure?  The answer lies with the stacky group
structure.

The space of functions 
on every $G$-$H$ bibundle
naturally acquires the structure of an $\calA(G)$-$\calA(H)$ bimodule.
Composition of bibundles corresponds to tensor product of bimodules,
so, under suitable technical assumptions, we get a functor from 
$\LieGrpdPrBibu$ to a category in which the objects are algebras and
the morphisms are bimodules.  
By analogy with the
usual Gelfand functor, we view it as contravariant. Since this functor takes
categorical products of groupoids to tensor products of algebras,
which are not categorical products, 
it certainly does not take group objects to group objects.  In fact,
it does not even take them to Hopf algebras, as we saw already for the
example of the noncommutative torus.  (There, the functor takes the 
unit of the group to a $\bbC$-$\calA$ bimodule rather than to a
homomorphism $\bbC \leftarrow \calA$.)  Instead,
the 
image of a stacky group is another structure which we call a
{\bf hopfish algebra}
\cite{ta-we-zh:hopfish,bl-ta-we:hopfish}.

A space of functions on the bibundle $E_m$  of multiplication for a
stacky group with function algebra $\calA$
becomes an $(\calA \otimes \calA)$-$\calA$ bimodule, which we denote by
$\Bimo{\Delta}$ and view as the bimodule of
comultiplication.  Analogously, the bibundle $E_e$ of the unit
element, already mentioned above, 
becomes the $\bbC$-$\calA$ bimodule $\Bimo{\epsilon}$ of the
counit. By functoriality these bibundles  satisfy coassociativity and
counitality relations
\begin{gather}
 (\calA \otimes \Bimo{\Delta} )
 \otimes_{\calA \otimes \calA} \Bimo{\Delta}
 \cong
 (\Bimo{\Delta} \otimes \calA )
 \otimes_{\calA \otimes \calA} \Bimo{\Delta}
 \,,\\
 (\Bimo{\epsilon} \otimes  \calA )
 \otimes_{\calA \otimes  \calA} \Bimo{\Delta}
 \cong \calA \cong
 (\calA \otimes \Bimo{\epsilon} )
 \otimes_{\calA \otimes  \calA} \Bimo{\Delta}
 \,,
\end{gather}
up to isomorphism of bimodules. We thus obtain a weak comonoidal
object in the weak 2-category of algebras and bimodules, which is also
called a sesquilinear sesquialgebra.

The algebraic image of the bibundle $E_{\inv}$ of the inverse  is more
difficult to interpret.
Because $E_{\inv}$ is a $G$-$G$ bibundle, the space of
functions on $E_{\inv}$ is an $\calA$-$\calA$ bimodule, which we would
expect to become the bimodule of an
antipode. But a Hopf antipode is an algebra antihomomorphism, so
the hopfish antipode should be an
$\calA$-$\calA^\op$ bimodule.  A way to accomplish this is to use the
star structure on $\calA$ which comes from the groupoid inverse of $G$,
$a^*(g) = \overline{a(g^{-1})}$, in order to convert the
$\calA$-$\calA$ bimodule of functions on $E_{\inv}$ into the
$\calA$-$\calA^\op$ bimodule $\Bimo{S}$ of the hopfish antipode.
Just as we illustrated above in the case of Poisson groups with
the diagram \ref{diagram-opinverse}, it is not easy to formulate what
is required of the antipode in a hopfish algebra.  
But a definition is given in \cite{ta-we-zh:hopfish}, with a modification
in \cite{bl-ta-we:hopfish}
to cover the case of  noncommutative torus algebras, and the bimodule
described above does satisfy the definition.  

Like an ordinary coproduct, the hopfish
coproduct bimodule $\Bimo{\Delta}$ 
can be used to multiply two right $\calA$-modules $T, T' \in
\mathrm{Mod}_{\calA}$ by
\begin{equation}
 T \otimes_{\Bimo{\Delta}} T' := (T \otimes T')
 \otimes_{\calA \otimes \calA} \Bimo{\Delta} \,.
\end{equation}
In \cite{bl-ta-we:hopfish},
we have tried out this new multiplication on certain modules over the 
hopfish algebra associated to the stacky group $S^1/\!/\bbZ$.
For $\alpha \in \bbR$, $p,q \in \bbZ$ relatively prime, we have shown
there that
\[
 T^{\alpha}_{pq} := \calA/(\E^{-\I\alpha} \Abas_{pq} - 1) \calA \,,
\]
is a simple module generated by $\xi := [1]$ with $a_{pq} \cdot \xi =
e^{i \alpha} \xi$, and that adding any integer multiple of $\lambda$ to
$\alpha$ results in an isomorphic module. 
(For $p$, $q$ not relatively prime the situation is
slightly more complicated.) Some calculations lead to the following
result.

\begin{Theorem}
\label{th:TensorProduct}For $p_1 \neq 0$ or $p_2 \neq 0$ we have:
\[
 T^{\alpha_1}_{p_1 q_1} \Motimes T^{\alpha_2}_{p_2 q_2} \cong
 \gcd(p_1,p_2)\, T^{\alpha}_{p q} \,,
\]
where
\begin{equation}
\label{eq:alphapqDef}
 p := \lcm(p_1,p_2) \,,~
 q := \frac{p_1 q_2 + p_2 q_1}{\gcd(p_1, p_2)} \,,~
 \alpha := \frac{\alpha_1 p_2 + \alpha_2 p_1}{\gcd(p_1,p_2)} \,.
\end{equation}
For $p_1 = 0$ and $p_2 = 0$ we have:
\[
 T^{\alpha_1}_{0, q_1} \Motimes T^{\alpha_2}_{0, q_2} \cong
 \begin{cases}
 T^{\alpha}_{0, q} \,, &\mathrm{for}\quad
 \frac{\alpha_1 q_2 - \alpha_2 q_1}{\lambda \gcd(q_1,q_2)} \in \bbZ \bmod \lcm(q_1,q_2)\frac{2\pi}{\lambda}\\
 0 \,,& \mathrm{otherwise}
 \end{cases} \,,
\]
where
\begin{equation}
\label{eq:alphapqDef2}
q := \gcd(q_1,q_2) \,, \quad
 \alpha := s_1 \alpha_2 - s_2 \alpha_1
 \,,\quad s_1,s_2 \in \bbZ \,:\quad
 \frac{s_1 q_2 - s_2 q_1}{\gcd(q_1,q_2)} = 1 \,.
\end{equation}
\end{Theorem}

Observe that $q/p = q_1/p_1+q_2/p_2$ and
$\alpha/p=\alpha_1/p_1+\alpha_1/p_2.$  Surprisingly, we did not
notice these simple relations until we found the following geometric
interpretation of them.

When we identify $\calA$ with an algebra of functions on the torus
$\bbT^2$, with coordinates $(\theta_1,\theta_2)$, the algebra element
$e^{-i\alpha} a_{pq}-1$ becomes the function
$e^{-i\alpha}e^{i(p\theta_1+q\theta_2)}-1$, and so the classical
  analogue of the quotient module $T^\alpha_{pq}=
\calA/(e^{-i\alpha} a_{pq}-1)\calA$
  appears to be the functions on the embedded circle
  in $\bbT^2$ defined by the equation $p\theta_1+q\theta_2-\alpha=0$,
  or $\theta_1=-\frac{q}{p}\theta_2+\frac{\alpha}{p}.$  The classical
  analogue of the hopfish structure on $\calA$ (and hence the
  symplectic model of the group structure on $S^1/\bbZ$) turns out to
  be the symplectic group{\em oid} structure $\bbT^2\rightrightarrows
  \bbT^1$
 for which the
  source and target maps are the projection in the $\theta_2$
  direction, and the composition law is addition in $\theta_1$.  The
  tensor product operation on modules corresponds to the application
  of the composition law on the groupoid to the embedded circles which
  represent them.  This results in the addition of fractions mentioned
  above.

Note that, if we ``unwrap'' the $\theta_2$ circle to a line, we obtain
the cotangent bundle $\bbT^1 \times \bbR \cong T^*\bbT^1,$ the
symplectic groupoid $\Gamma(\bbT^1)$ of $\bbT^1$ with the zero Poisson
structure.  The groupoid structure described in the preceding
paragraph is the second symplectic groupoid structure on $\Gamma(\bbT^1)$
which is obtained by lifting the (Poisson) Lie group structure on
$\bbT^1$ given by addition in $\theta$.  This is an instance of the
double symplectic groupoid structures attached to general Poisson Lie
groups  \cite{lu-we:groupoides}.

What has become of the parameter $\lambda$?  In fact, it seems that
the curve $p\theta_1+q\theta_2-\alpha$ represents only one generator
of the module in question.  The others are obtained by applying the
unitary basis elements of the algebra itself.  These correspond to
constant ``integer'' bisections of $\Gamma(\bbT^2)\cong T^*\bbT^2$;
the bisection $m d\theta_1 + n d\theta_2$ 
acts via the source and target maps, determined
by the Poisson structure (see \cite{we:rotation}), as translation by
$(-\lambda n, \lambda m)$.  This yields the collection of all circles
of the form $p\theta_1+q\theta_2-(\alpha+\lambda(np-mq))=0$.  
Under the assumption $\gcd(p,q)=1$, all the integer multiples of
$\lambda$ occur, so that this collection of circles, like the
isomorphism class of $T^\alpha_{pq}$, depends only on the image of
$\alpha$ in $S^1/\bbZ$. 

%

\end{document}